\documentclass[12pt]{article}

\usepackage{amssymb}

\newcommand{\bea}{\begin{eqnarray}}
\newcommand{\eea}{\end{eqnarray}}
\newcommand{\beay}{\begin{eqnarray*}}
\newcommand{\eeay}{\end{eqnarray*}}

\newcommand{\nlim}{\lim_{n\rightarrow\infty}\ }
\newcommand{\nlinf}{\liminf_{n\rightarrow\infty}\ }

\newcommand{\zlim}{\lim_{t\rightarrow 0}\ }

\newcommand{\gosto}{\rightarrow}
\newcommand{\sgn}{\mbox{sgn}\,}
\newcommand{\sech}{\mbox{sech}}
\newcommand{\proof}{\noindent{\bf Proof}.\quad}

\newcommand{\eq}[1]{(\ref{eq#1})}
\newcommand{\thurem}[1]{\ref{thurem#1}}

\newcommand{\prop}[1]{\ref{prop#1}}
\newcommand{\conjecture}[1]{\ref{conjecture#1}}
\newcommand{\keywords}{\noindent{\bf Key words.}\ }
\newcommand{\subclass}{\noindent{\bf AMS subject classifications.}\ }

\newtheorem{theorem}{Theorem}[section]
\newtheorem{lemma}[theorem]{Lemma}
\newtheorem{corollary}[theorem]{Corollary}

\newtheorem{example}[theorem]{Example}

\newtheorem{conj}[theorem]{Conjecture}
\newtheorem{proposition}[theorem]{Proposition}

\newtheorem{remark}[theorem]{Remark}

\begin{document}

\title{On the global attractivity and oscillations in a class of second order difference equations from macroeconomics}         

\author{ Hassan A. El-Morshedy\thanks{Department of Mathematics, Damietta Faculty of Science, New Damietta 34517, Egypt.} \\({\tt elmorshedy@yahoo.com})}        

\date{}
\maketitle

\begin{abstract}
New global attractivity criteria are obtained for  the second order difference equation
\[
x_{n+1}=cx_{n}+f(x_{n}-x_{n-1}),\quad n=1,\,2,\, \ldots
\]
via a Lyapunov-like method. Some of these results are sharp and support recent related conjectures.
Also, a necessary and sufficient condition for the oscillation of this equation is obtained using comparison with a second order linear difference equation with positive coefficients.
\end{abstract}

\keywords  Nonlinear difference equations, Global attractivity, Oscillation, Macroeconomics models

\subclass 39A10, 39A11

\section{Introduction}
Consider the second order difference equation
\begin{equation}\label{eq1}
x_{n+1}=cx_{n}+f(x_{n}-x_{n-1}),\quad n=1,\,2,\, \ldots
\end{equation}
where $c\in[0,\,1)$, $f:\mathcal{R}\gosto \mathcal{R}$ is a continuous real function and the initial values $x_{0},\,x_{1}$ are real numbers. Various particular cases of \eq{1} have appeared in mathematical models of macroeconomics. For prototype examples, the reader is referred to Samuelson \cite{sam}, Hicks \cite{hicks} and Puu \cite{puu}. Motivated by those examples Sedaghat in \cite{s1997} proposed and investigated the general form \eq{1}. We mention here that for sigmoidal or tanh-like nonlinearities, equation \eq{1} can also be regarded as the discrete analogue of the single delayed neuron model
\[
x'(t)=-\alpha x(t)+f(x(t)-x(t-\tau))
\]
using (forward) Newton discretization scheme with step size equals $\tau$. An account of the stability analysis and/or the oscillations of the above continuous {\it neuronic} equation and some related equations can be found in \cite{elmoel,elmog,gl1,gh} while a higher order discrete neuronic version has been investigated by \cite{hamaya}.

The global attractivity (stability), boundedness and/or oscillations of \eq{1} have been considered by \cite{kents,s1997,s2002,s2005}. Very recently, Li and Zhang \cite{lizhang} studied its bifurcation .

It will be assumed, without loss of generality, that the origin is the unique equilibrium point of \eq{1}.

The known global attractivity results for \eq{1} are collected in the following result.
\begin{theorem}\label{thurem1}
Assume that $|f(t)|\leq a|t|$ for all $t$. The origin of \eq{1} is globally attracting if any one of the following conditions is satisfied:
\begin{enumerate}

\item [(c1)]\cite[Corollary 1]{kents}
    \begin{equation}\label{eq2}
 a<\frac{1-c}{2}.
 \end{equation}
\item [(c2)] \cite[Corollary 4]{kents} $tf(t)\geq 0$ for all $t\in \mathcal{R}$ and
  \begin{equation}\label{eq4}
	a<\max\{b,\,1-c,\,d\}
   \end{equation}
where $b=(1-\sqrt{1-c})^{2}$ and $d=\frac{2-c}{3-c}$. The same conclusion also holds when $a=1-c$, $c\not=0$ or $a=d$.
  \item [(c3)] \cite[Theorem 8]{s2005} $f(t)\geq 0$ for all $t\in \mathcal{R}$ and
	\begin{equation}\label{eq5}
	a<\max\{1-c,\,c\}.
   \end{equation}
\end{enumerate}
\end{theorem}

\begin{conj}\label{conjecture1}
\cite{kents,s2005} If $|f(t)|\leq a|t|$ and $tf(t)\geq 0$ for all $t\in \mathcal{R}$ where $a\in(0,\,1)$, then the origin of \eq{1} is globally attracting.
\end{conj}

\begin{conj}\label{conjecture2}
\cite{s2005} If $0\leq f(t)\leq a|t|$ for some $a\in(0,\,1)$ and all $t\in \mathcal{R}$, then the origin of \eq{1} is globally attracting.
\end{conj}

In this work, we contribute to the validity of the above conjectures by improving Theorem \thurem{1}. We use a Lyapunov-like method to investigate the global attractivity of the origin of \eq{1}. Moreover, the oscillation of \eq{1} is studied via comparison with the oscillation of a second order difference equation with constant coefficients. This helps us to improve \cite[Theorem 2(a)]{kents}. Necessary and sufficient condition for the oscillation of \eq{1} is obtained when $\frac{f(x)}{x}$ is maximized at $0$. Here  an equation is called oscillatory if each of its solutions is neither eventually negative nor eventually positive. 

The following proposition will be needed in  some of the proofs below.
\begin{proposition}\label{prop1}
Let $\{x_{n}\}$ and $\{y_{n}\}$ be two real sequences such that $y_{n}=x_{n}-x_{n-1}$, $n\geq 1$ and $x_{2}=c x_{1}+f(x_{1}-x_{0})$. Then $\{x_{n}\}$ is a solution of equation \eq{1} if and only if $\{y_{n}\}$ is a solution of the equation
\begin{equation}\label{eq6}
y_{n+1}=cy_{n}+f(y_{n})-f(y_{n-1}), \quad n\geq 2.
\end{equation}
Moreover, the origin of \eq{1} is globally attracting if and only if the origin of \eq{6} is globally attracting.
\end{proposition}
\proof If $\{x_{n}\}$ is a solution of \eq{1} and $y_{n}=x_{n}-x_{n-1}$ for $n\geq 1$, it follows that
\[
y_{n+1}=x_{n+1}-x_{n}=cy_{n}+f(y_{n})-f(y_{n-1}), \quad n\geq 2.
\]
Now assume that $y_{n}$ is a solution of \eq{6}, then
\[
x_{n+1}-cx_{n}-f(x_{n}-x_{n-1})=x_{n}-cx_{n-1}-f(x_{n-1}-x_{n-2}),\quad n\geq 2.
\]
Since $x_{2}-c x_{1}-f(x_{1}-x_{0})=0$, then the above equality implies that $\{x_{n}\}$ is a solution of \eq{1} as desired. This proves the first part of the proposition. For the second part, we prove only that the origin of \eq{1} is globally attracting provided that the origin of \eq{6} is globally attracting. This clearly follows from \eq{1} since for any solution $\{x_{n}\}$ of \eq{1} there exists a solution $\{y_{n}\}$ of \eq{6} such that
\[
x_{n+1}=cx_{n}+f(y_{n}),\quad n\geq 1
\]
which implies that $y_{n+1}-f(y_{n})=(c-1)x_{n}$ and hence 
\[
\nlim x_{n} = \nlim (c-1)^{-1}\{y_{n+1}-f(y_{n})\}=0.
\]

\setcounter{equation}{0}
\section{Global attractivity}We start with the following sharp result.

\begin{theorem}\label{thurem2}
Assume that $|f(t)|\leq a|t|$ for all $t$. If
\begin{equation}\label{eq7}
a<\frac{1+c}{2},
\end{equation}
then the origin of \eq{1} is globally attracting.
\end{theorem}

\proof
Due to Proposition \prop{1}, it is enough to prove the global attractivity of the origin of \eq{6}.

Let $\{V_{n}\}_{n\geq 2}$ be defined as follows
\[
V_{n}=\beta (y_{n}-f(y_{n-1}))^{2}+\gamma (f(y_{n-1}))^{2},\quad n\geq 2
\]
where $\{y_{n}\}$ be any solution of \eq{6} and $\beta,\,\gamma$ are positive real numbers to be determined later. Then
\bea\label{eq8}
\Delta V_{n}&=&\beta (c^{2}-1)y_{n}^{2}-2 \beta (c-1)y_{n}f(y_{n-1})+\gamma (f(y_{n}))^{2}-\gamma (f(y_{n-1}))^{2}\nonumber\\
&\leq & [\gamma a^{2}-\beta (1-c^{2})]y_{n}^{2}+2 \beta (1-c)y_{n}f(y_{n-1})-\gamma (f(y_{n-1}))^{2}, \ n\geq 2.
\eea
Completing square with respect to $f(y_{n-1})$, it follows that
\bea\label{eq9}
\Delta V_{n} &\leq & -Ay_{n}^{2}-\gamma (f(y_{n-1})-\frac{\beta}{\gamma}(1-c)y_{n})^{2}\nonumber\\
&\leq & -A y_{n}^{2},\quad n\geq 2
\eea
where $A=-\gamma a^{2}-\frac{\beta^{2}}{\gamma}(1-c)^{2}+\beta (1-c^{2})$. We require that $A>0$ for some $\gamma,\,\beta>0$. This is equivalent to saying that
\[
\frac{1+c-\sqrt{(1+c)^{2}-4a^{2}}}{2(1-c)}<\frac{\beta}{\gamma}<\frac{1+c+\sqrt{(1+c)^{2}-4a^{2}}}{2(1-c)}.
\]
So $\beta$ and $\gamma$ exist if \eq{7} holds. Hence summing \eq{9} from $2$ to $n$, we obtain
\[
V_{n+1}-V_{2}\leq -A\sum_{i=2}^{n}y_{i}^{2}.
\]
Since $V_{n}\geq 0$ for all $n\geq 2$, the above inequality leads to the convergence of  $\sum_{i=2}^{\infty}y_{i}^{2}$  and hence $\nlim y_{n}=0$.

\begin{remark}\label{rem1}
Generally; the above result can not be weakened for functions satisfying $|f(t)|\leq a|t|$ on $\mathcal{R}$. Indeed when $a\geq \frac{1+c}{2}$, one can find functions $f$ with which equation \eq{1} has solutions that are not attracted to the origin. For example, when $f(t)=-at$ (see \cite[p.1261]{kents}), the characteristic polynomial associated with equation \eq{1} has negative solution $\lambda\leq -1$ and so for $x_{n}=\lambda^{n}$, the solution $\{x_{n}\}$ diverges. Moreover, we observe that Theorem \thurem{2} improves Theorem \thurem{1}(c1) for $c\in (0,\,1)$ and  Theorem \thurem{1}(c2),(c3) when $c>\frac{1}{3}$.
\end{remark}
\begin{theorem}\label{thurem3}
Assume that $0\leq f(t)\leq a|t|$ for all $t\in \mathcal{R}$. If
\begin{equation}\label{eq10}
a^{2}<\frac{1+c}{2l}, \quad l=\max\{\frac{a}{a+c},\,\frac{a}{a+1-c}\},
\end{equation}
then the origin of \eq{1} is globally attracting.
\end{theorem}
\proof Since
\[
x_{n+1}=cx_{n}+f(x_{n}-x_{n-1})\geq cx_{n},\quad n\geq 2,
\]
then $\{x_{n}\}$ is eventually of one sign. Let $x_{n}<0$ for all $n\geq n_{0}>2$. Then $0>x_{n+1}> x_{n}$ for all $n\geq n_{0}$ and hence $\nlim x_{n}=0$. Therefore, $x_{n}>0$ for all $n\geq n_{0}$ which implies that
\[
x_{n+1}<x_{n}+f(x_{n}-x_{n-1}), \quad n\geq n_{0}+1.
\]
That is,
\begin{equation}\label{eq11}
y_{n+1}<f(y_{n}),\ \mbox{for all }n\geq n_{0}+1
\end{equation}
where $y_{n}=x_{n}-x_{n-1}$ for $n\geq n_{0}+1$. This inequality yields
\[
y_{n}<a|y_{n-1}|, \quad n\geq n_{0}+2
\]
or equivalently
\bea\label{eq12}
\frac{y_{n}}{a}<\left\{\begin{array}{ll}
y_{n-1},\quad & y_{n-1}>0\\
-y_{n-1},\quad & y_{n-1}<0.
\end{array}
\right.
\eea
Using \eq{6} and \eq{11}, we obtain
\beay
y_{n}&=&cy_{n-1}+f(y_{n-1})-f(y_{n-2})\\
&<&\left\{\begin{array}{ll}
cy_{n-1}+f(y_{n-1}),\quad & y_{n-1}<0\\
(c-1)y_{n-1}+f(y_{n-1}),\quad & y_{n-1}>0.
\end{array}
\right.
\eeay
So \eq{12} yields
\beay
y_{n}<\left\{\begin{array}{ll}
\frac{-c}{a}y_{n}+f(y_{n-1}),\quad & y_{n-1}<0\\
\frac{c-1}{a}y_{n}+f(y_{n-1}),\quad & y_{n-1}>0.
\end{array}
\right.
\eeay
Rearranging,
\bea\label{eq13}
y_{n}&<&\left\{\begin{array}{ll}
\frac{a}{a+c}f(y_{n-1}),\quad & y_{n-1}<0\\
\frac{a}{a-c+1}f(y_{n-1}),\quad & y_{n-1}>0
\end{array}
\right.\nonumber\\
&<&l f(y_{n-1}), \quad n\geq n_{0}+2.
\eea
Define $V_{n}$ as in the proof of Theorem \thurem{2}, then \eq{8} and \eq{13} imply that
\[
\Delta V_{n}< [\gamma a^{2}-\beta (1-c^{2})]y_{n}^{2}+[2\beta (1-c) l - \gamma](f(y_{n-1}))^{2}.
\]
In view of \eq{10}, we have $\frac{1-c^{2}}{a^{2}}>2(1-c)l$. Thus the values of $\beta,\,\gamma$ can be chosen such that
\[
2(1-c)l<\frac{\gamma}{\beta}<\frac{1-c^{2}}{a^{2}}
\]
which yields
\[
\gamma a^{2}-\beta (1-c^{2})=-B<0\quad \mbox{and}\quad 2\beta (1-c) l - \gamma<0.
\]
Therefore,
\[
\Delta V_{n}<-B y_{n}^{2}, \quad n\geq n_{0}+2
\]
which implies, as in the proof of Theorem \thurem{2}, that $\nlim x_{n}=0$.

\begin{remark}
It is easy to see that \eq{10} improves \eq{7} for $a\in (0,\,1]$. Moreover, since
\[
\frac{1+c}{2l}=\left\{\begin{array}{ll}
\frac{(1+c)(a+c)}{2a},\quad & c\leq \frac{1}{2}\\
\frac{(1+c)(1+a-c)}{2a},\quad & c\geq \frac{1}{2}
\end{array}
\right.\nonumber\\
\]
Then condition \eq{10} holds, with $a\in (0,\,1]$,  provided that 
\[
a^{2}<\left\{\begin{array}{ll}
\frac{(1+c)^{2}}{2},\quad & c\leq \frac{1}{2}\\
\frac{(1+c)(2-c)}{2},\quad & c\geq \frac{1}{2}
\end{array}
\right.\nonumber\\
\]
which in turn is satisfied for all $c\geq -1+\sqrt{2}$. This gives a partial positive answer to Conjecture  \conjecture{2}.

On the other hand, Theorem \thurem{3} can be used when $a>1$. For example when $c=1/2$, \eq{10} becomes $a^{3}-\frac{3}{4}a-\frac{3}{8}<0$ which is satisfied for all $a<1.0519...$. This suggests that Conjecture \conjecture{2} is true for $a\in(0,\,\alpha)$ and some $\alpha>1$.
\end{remark}
\setcounter{equation}{0}
\section{Oscillations}
The following result refines Lemma 3 in \cite{kents} for certain types of $f$.
\begin{lemma}\label{lem1}
Assume that $tf(t)\geq 0$ for all $t\in \mathcal{R}$ and $\{x_{n}\}$ be any nonoscillatory solution of equation \eq{1}. If
\begin{equation}\label{eq14}
|f(t)|\leq a|t| \quad \mbox{for all } |t|\geq t_{0}\ \mbox{and some } a\leq 1
\end{equation}
where $t_{0}$ is a sufficiently large number, then $x_{n}\Delta x_{n}<0$ eventually.
\end{lemma}
\proof Using \cite[Lemma 3]{kents}, we see that $\{x_{n}\}$ is eventually monotonic. We assume that $x_{n}>0$ for all $n\geq n_{0}>0$ (the case when $x_{n}<0$, eventually, can be handled similarly).

For the sake of contradiction, we assume that $\Delta x_{n}>0$, $n\geq n_{1}\geq n_{0}$. It follows that either $\nlim x_{n}=l>0$ or $\nlim x_{n}=\infty$. The first case is impossible as the only possible limit of $\{x_{n}\}$ is zero. Now, in view of the increasing nature of $\{x_{n}\}_{n\geq n_{1}}$, equation \eq{1} implies that
\[
(1-c) x_{n+1}<f(x_{n}-x_{n-1}),\quad n\geq n_{1}.
\]
So $\nlim f(\Delta x_{n-1})=\infty$ which is possible only if $\nlim \Delta x_{n-1}=\infty$. Using \eq{14}, it is easy to find $n_{2}\geq n_{1}$ such that $f(\Delta x_{n-1})\leq a\Delta x_{n-1}$ for $n\geq n_{2}$. Therefore, equation \eq{1} yields
\beay
x_{n+1}&\leq& cx_{n}+a (x_{n}-x_{n-1})\\
&<& x_{n}+a \Delta x_{n-1},\quad n\geq n_{2}.
\eeay
Thus $\Delta x_{n} < a \Delta x_{n-1}\leq \Delta x_{n-1}$ for $n\geq n_{2}$ and hence $\nlim \Delta x_{n} \not=\infty$.

\begin{remark}\label{remark2}
We observe that condition \eq{14} covers many types of functions. For example; each of the following functions satisfies \eq{14}:
\[
f:\ |f(t)|\leq a \quad \mbox{for all } t \ \mbox{ and some } a>0,
\]
\begin{equation}\label{eq15}
f:\ |f(t)|\leq a|t|\quad \mbox{for all } t \ \mbox{ and some } a\leq 1,
\end{equation}
and the sublinear function
\begin{equation}\label{eq16}
f: \ f(t)=|t|^{\lambda} \sgn t,\quad \lambda \in (0,\,1)
\end{equation}
where
\beay
\sgn t=\left\{\begin{array}{ll}
-1,\quad & t<0\\
0, \quad &t=0\\
+1, \quad & t>0
\end{array}
\right.
\eeay
\end{remark}

\begin{theorem}\label{thurem4}
Assume that all assumptions of Lemma~\ref{lem1} hold and there exist $\alpha_{1},\,\alpha{_2}\in(0,\,\infty)$ such that $\liminf_{t\rightarrow 0^{-}}\ \frac{f(t)}{t}\geq \alpha_{1}$ and $\liminf_{t\rightarrow 0^{+}}\ \frac{f(t)}{t}\geq \alpha_{2}$. If
\begin{equation}\label{eq17}
(1-\sqrt{1-c})^{2}<\alpha_{i}<(1+\sqrt{1-c})^{2},\quad \mbox{for}\ i=1,2
\end{equation}
then equation \eq{1} is oscillatory.
\end{theorem}
\proof If $\{x_{n}\}$ is a nonoscillatory solution of \eq{1}, then it is either eventually negative or eventually positive. Assume that $\{x_{n}\}$ is  eventually positive. Then $x_{n}>0$ for all $n\geq \bar{n}$ for some $\bar{n}>1$ and Lemma~\ref{lem1} implies that $\Delta x_{n}<0$ for all $n\geq n_{1}>\bar{n}$ which in turn yields $\nlim x_{n}=0$, $\nlim \Delta x_{n-1}=0$ and consequently $\nlinf \frac{f(\Delta x_{n-1})}{\Delta x_{n-1}} \geq \alpha_{1}$. Therefore, for any $\epsilon>0\,(\alpha_{1}-\epsilon>0)$ there exists $n_{2}\geq n_{1}$ such that
\[
f(\Delta x_{n-1})\leq (\alpha_{1}-\epsilon)\Delta x_{n-1}, \quad n\geq n_{2}.
\]
Substituting into \eq{1},
\[
x_{n+1}\leq (c+\alpha_{1}-\epsilon)x_{n}-(\alpha_{1}-\epsilon)x_{n-1},
\]
which can be rewritten in the self-adjoint form
\[
\Delta ((\alpha_{1}-\epsilon)^{-n+1}\Delta x_{n-1})+(1-c)(\alpha_{1}-\epsilon)^{-n}x_{n}\leq 0, \quad n\geq n_{2}.
\]
The existence of a positive solution of the above inequality implies (see \cite[p.470]{ce}) the nonoscillation of the second order difference equation
\[
\Delta ((\alpha_{1}-\epsilon)^{-n+1}\Delta z_{n-1})+(1-c)(\alpha_{1}-\epsilon)^{-n}z_{n}= 0, \quad n\geq n_{2}.
\]
This is possible only  if the corresponding characteristic polynomial
\[
\lambda^{2}-(c+\alpha_{1}-\epsilon)\lambda +(\alpha_{1}-\epsilon)=0
\]
has positive solutions; that is when $(\alpha_{1}-\epsilon)\leq (1-\sqrt{1-c})^{2}$ or $(\alpha_{1}-\epsilon)\geq (1+\sqrt{1-c})^{2}$. Since $\epsilon$ is arbitrary, then these inequalities can not hold due to \eq{17}. This contradiction proves this case.

When $x_{n}<0$ eventually, similar arguments lead to the proof. We omit the details to avoid repetition.

\begin{remark}\label{remark3}
Since $\alpha_{1}$ and $\alpha_{2}$ are calculated at zero, then Theorem \thurem{4} improves \cite[Theorem 2(a)]{kents}.
\end{remark}

Theorem 2(b) in \cite{kents} asserts that \eq{1} is nonoscillatory if $tf(t)\geq 0$, $|f(t)|\leq a|t|$ for all $t\in \mathcal{R}$ and $a\leq b=(1-\sqrt{1-c})^{2}$. Now, if $f$ is continuously differentiable at zero and $\frac{f(x)}{x}$ is maximized at zero (i.e., $\frac{f(x)}{x}\leq f'(0)$ for all $x\not=0$), then $\alpha_{1}=\alpha_{2}=f'(0)$ and a combination of Theorem \thurem{4} and Theorem 2(b) in \cite{kents} leads to the following necessary and sufficient condition for the oscillation of \eq{1}.

\begin{corollary}\label{corollary1}
Assume that $f$ is continuously differentiable at zero and $0<\frac{f(t)}{t}\leq f'(0)<1$ for $t\not=0$. Then equation \eq{1} is oscillatory if and only if $f'(0)>(1-\sqrt{1-c})^{2}$.
\end{corollary}

\begin{example}\label{exam1}
Consider the discrete single neuron model
\[
x_{n+1}=cx_{n}+a \tanh (x_{n}-x_{n-1}), \quad 0<a<1.
\]
Here $f(t)=a \tanh t$ and $f'(t)\leq a \sech^{2}0 =a$ for all $t$. Then all solutions of this model oscillate if and only if $a>(1-\sqrt{1-c})^{2}$.
\end{example}

Theorem \thurem{4} can also be used to investigate the oscillation of \eq{1} when $f$ is a sublinear function. In this case, in addition to \eq{14}, $f$ satisfies
\begin{equation}\label{eq18}
\zlim \frac{f(t)}{t}=\infty.
\end{equation}
Therefore, $\alpha_{1},\,\alpha_{2}$ can be suitably chosen to satisfy \eq{17}.
\begin{corollary}\label{corollary2}
Assume that $f$ satisfies \eq{18} and all assumptions of Lemma~\ref{lem1}. Then equation \eq{1} is oscillatory.
\end{corollary}

{\bf Acknowledgments.} The author would like to thank  the referees  for their valuable comments and suggestions.



\begin{thebibliography}{10}

\bibitem{ce}
S. Chen and Lynn H. Erbe, {\em Riccati techniques and discrete oscillations}, J. Math. Anal. Appl. 142(1989), 468-487. 

\bibitem{elmoel}
 H. A. El-Morshedy and B. M. El-Matary, { \em Oscillation and global asymptotic stability of a
neuronic equation with two delays}, Electron. J. Qual. Theory Differ. Equ. 2008(2008), 1-21.

\bibitem {elmog}
H. A. El-Morshedy and K. Gopalsamy, {\em On the oscillation and asymptotic behaviour of solutions of a neuronic equation},  Funkcial. Ekvac. 44(2001) 83-98.

\bibitem {gl1}
 K. Gopalsamy and I. K. C. Leung, {\em Convergence under dynamical thresholds with delays}, IEEE Trans. Neural Networks, 8(1997),  341-348.

\bibitem{gh}
I. Gy\"{o}ri and F. Hartung, {\em Stability of a single neuron model with delay}, J. Comput. Appl. Math. 157(2003), 73-92.

\bibitem{hamaya}
Y. Hamaya, {\em On the asymptotic behaviour of solutions of neuronic difference equations}, Proc. of Internat. Conf. on Difference Eqns, Special Funcs. and Appl. (S. Elaydi et al, Editors) Munich, Germany, 25-30 July 2005, World Scientific (2005), pp. 258-265.

\bibitem{hicks}
J. R. Hicks, {\em A contribution to the theory of the trade cycle}, 2nd Ed., Clarendon Press, Oxford, 1965.

\bibitem{kents}
C. M. Kent and H. Sedaghat, {\em Global stability and boundedness in $x_{n+1}=cx_{n}+f(x_{n}-x_{n-1})$}, J. Difference Equ. Appl. 10(2004), 1215-1227.

\bibitem{lizhang}
S. Li and W. Zhang, {\em Bifurcation in a second-order difference equation from macroeconomics}, J. Difference Equ. Appl. 14(2008), 91-104.

\bibitem{puu}
T. Puu, {\em Nonlinear economic dynamics}, 3rd Ed., Springer, New York, 1993.

\bibitem{sam}
H. A. Samuelson, {\em Interaction between the multiplier analysis and the principle of acceleration}, Rev. Econ. Stat. 21(1939), 75-78.

\bibitem{s1997}
H. Sedaghat, {\em A class of nonlinear second order difference equations from macroeconomics}, Nonlinear Anal. 29(1997), 593-603.


\bibitem{s2002}
H. Sedaghat, {\em Regarding the equation $x_{n+1}=cx_{n}+f(x_{n}-x_{n-1})$}, J. Difference Equ. Appl. 8(2002), 667-671.


\bibitem{s2005}
H. Sedaghat, {\em Global attractivity, oscillation and chaos in a class of nonlinear second order difference equations}, Cubo 7(2005), 89-110.




\end{thebibliography}
\end{document}